\pdfoutput=1

\documentclass[12pt]{article}\textwidth 160mm\textheight 235mm
\usepackage{amsfonts}
\usepackage{amssymb}
\usepackage{graphicx}
\usepackage{amsmath}
\usepackage{tikz}
\usepackage{cancel}
\usepackage{todonotes}
\usetikzlibrary{matrix}
\usetikzlibrary{arrows,shapes}

\def\ve{\varepsilon}
\def\epsilon{\varepsilon}

\def\emp{\emptyset}

\def\dom{{\rm dom}\,}

\def\O{{\cal O}}

\def\ox{\overline{x}}

\def\disp{\displaystyle}

\def\tto{\;{\lower 1pt \hbox{$\rightarrow$}}\kern-10pt
\hbox{\raise 2pt\hbox{$\rightarrow$}}\;}
\def\Hat{\widehat}

\def\Bar{\overline}
\def\ra{\rangle}
\def\la{\langle}
\def\ve{\varepsilon}

\def\h{\hfill\Box}
\def\R{\mathbb{R}}

\def\ox{\bar{x}}

\def\int{\mbox{\rm int}}

\def\dom{\mbox{\rm dom}\,}

\def\h{\hfill\triangle}
\def\dn{\downarrow}
\def\O{\Omega}
\def\ph{\varphi}
\def\emp{\emptyset}
\def\st{\stackrel}
\def\oR{\Bar{\R}}
\def\lm{\lambda}

\def\gg{\gamma}

\def\th{\theta}
\def\vt{\vartheta}
\def\sce{\setcounter{equation}{0}}

\begin{document}
\vspace*{0.5in}
\begin{center}
{\bf VARIATIONAL ANALYSIS AND VARIATIONAL RATIONALITY\\IN BEHAVIORAL SCIENCES: STATIONARY TRAPS}\\[2ex]
BORIS S. MORDUKHOVICH\footnote{Department of Mathematics, Wayne State University, Detroit, MI 48202, USA and RUDN University, Moscow 117198, Russia (boris@math.wayne.edu). Research of this author was partly supported by the National Science Foundation under grant DMS-1512846, by the Air Force Office of Scientific Research under grant \#15RT0462, and by the RUDN University Program 5-100.} and ANTOINE SOUBEYRAN\footnote{Aix-Marseille University (Aix-Marseille School of Economics), CNRS \& EHESS, Marseille 13002, France (antoine.soubeyran@gmail.com).}
\end{center}
\vspace*{0.05in}
\small{\bf Abstract.} This paper concerns applications of variational analysis to some local aspects of behavioral science modeling by developing an effective variational rationality approach to these and related issues. Our main attention is paid to local stationary traps, which reflect such local equilibrium and the like positions in behavioral science models that are not worthwhile to quit. We establish constructive linear optimistic evaluations of local stationary traps by using generalized differential tools of variational analysis that involve subgradients and normals for nonsmooth and nonconvex objects as well as variational and extremal principles.\\[1ex]

\newtheorem{Theorem}{Theorem}[section]
\newtheorem{Proposition}[Theorem]{Proposition}
\newtheorem{Remark}[Theorem]{Remark}
\newtheorem{Lemma}[Theorem]{Lemma}
\newtheorem{Corollary}[Theorem]{Corollary}
\newtheorem{Definition}[Theorem]{Definition}
\newtheorem{Example}[Theorem]{Example}
\renewcommand{\theequation}{{\thesection}.\arabic{equation}}
\renewcommand{\thefootnote}{\fnsymbol{footnote}}

\normalsize

\section{Introduction}\sce

Recent years have witnessed broad applications of advanced tools of variational analysis, generalized differentiation, and multiobjective (vector and set-valued) optimization to real-life models, particularly those related to economics and finance; see, e.g., \cite{aky18,bao-mor10,eic14,ktz15,m06,m18} and the references therein. Lately \cite{bms15a,bms15b}, certain variational principles and techniques have been developed and applied to models of behavioral sciences that mainly concern human behavior. The latter applications are based on the {\em variational rationality} approach to behavioral sciences initiated by Soubeyran in \cite{s09,s10}. Major concepts of the variational rationality approach include the notions of (stationary and variational) {\em traps}, which describe underlying positions of individual or group behavior related to making worthwhile decisions on changing or staying at the current position; see Section~2 for more details. Mathematically these notions correspond to points of equilibria, optima, aspirations, etc. and thus call to employ and develop powerful machinery of variational analysis and optimization theory for their comprehensive study and applications.

Papers \cite{bms15a,bms15b} mostly dealt with the study of {\em global variational traps} in connections with dynamical aspects of the variational rationality approach and its applications to goal systems in psychology \cite{bms15a} and to capability theory of wellbeing in behavioral sciences \cite{bms15b}. Appropriated tools of variational analysis developed in \cite{bms15a,bms15b} for the study and applications of the aforementioned global dynamical issues were related to set-valued extensions of the Ekeland variational principle in the case of set-valued mappings defined on quasimetric spaces and taken values in vector spaces with variable ordering structures; see also \cite{aky18,eic14,ktz15} for multiobjective optimization problems of this type.

The main goal of this paper is largely different from those in \cite{bms15a,bms15b}. Here we primarily focus on the study and variational analysis descriptions of {\em local stationary traps} in the variational rationality framework. It will be shown that the well-understood subdifferential notions for convex and nonconvex extended-real-valued functions as well as generalized normals to locally closed sets occur to be very instrumental to describe and characterize various types of such traps. Furthermore, the fundamental {\em variational} and {\em extremal principles} of variational analysis provide deep insights into efficient descriptions and evaluations of stationary traps and their approximate counterparts.\vspace*{0.05in}

The rest of the paper is organized as follows. In Section~2 we review relevant aspects of the variational rationality approach to behavioral sciences and formulate the major problems of our interest in this paper. Section~3 is devoted to global and local {\em evaluation aspects} of variational rationality. Here we introduce and investigate new notions of {\em linear evaluation} (both {\em optimistic} and {\em pessimistic}) of advantages and disadvantages to changing payoff and utility functions in various model structures. Then we discuss relationships between optimistic (resp.\ pessimistic) evaluations with {\em majorization-minimization} (resp.\ {\em minorization-maximization} algorithms of optimization. The main results here establish {\em complete subdifferential descriptions} of the {\em rates of change} for linear optimistic evaluations that are global for convex functions and (generally) local for nonconvex ones.

Section~4 conducts the study of {\em exact stationary traps} in behavioral dynamics. We discuss the meanings of such traps (which relate to a worthwhile choice between change and stay at the current position; see \cite{s18a,s18b} for more details) with considering both global and local ones while paying the main attention to local aspects. The usage of appropriate subgradients of variational analysis and linear optimistic evaluations allow us to derive efficient {\em certificates} for such stationary traps in general settings. Based on them, we reveal and discuss here remarkable special cases of such traps in behavioral models.

The next Section~5 deals with {\em approximate stationary traps}, which naturally appear in the framework of variational rationality and allow us to employ to their study powerful variational principles of variational analysis. Proceeding in this way, we derive here efficient descriptions of approximate stationary traps and establish verifiable conditions for the existence of such traps satisfying certain optimality and subdifferential properties.

Section~6 is devoted to {\em geometric} aspects of stationary trap evaluations by using generalized normals to nonconvex sets and their variational descriptions as well as by employing the fundamental {\em extremal principle} of variational analysis. Besides applications to behavior science modeling via the variational rationality approach, we provide here {\em behavioral interpretations} of generalized normals and the basic $\ve$-extremal principle around stationary trap positions. Furthermore, the important notion of {\em tilt perturbations} is involved in this section to evaluate stationary traps in behavioral dynamics associated with linear utility functions and the corresponding proximal payoffs.

The concluding Section~7 summaries our major finding in this paper. Then we discuss some open problems and directions of future research in this exciting area of interrelationships between variational analysis and variational rationality in behavioral sciences.

\section{Variational Rationality in Behavioral Sciences}\sce

The variational rationality approach to human behavior addresses a complex changing world under a long list of human limitations including limited cognitive resources, limited perception, unknown environment, etc.; see \cite{s09,s10} and the references therein for more details. All of this leads to the imperfect evaluation of payoffs. It has been well realized in behavioral sciences (see, e.g., \cite{kt79,tk91}) that agents evaluate their payoffs relative to moving reference points while depending on agents' experience and current environment. This definitely complicates the evaluation process and calls for possible simplifications.

In this paper we mostly concentrate on {\em local linear optimistic evaluations} (see below for the exact definitions) of nonlinear payoff and utility functions in the framework of the variational rationality approach to behavioral sciences. It will be done by using the basic tools and results of variational analysis and generalized differentiation.

First we briefly overview relevant aspects of variational rationality in the framework of {\em linear normed spaces}; this setting is assumed in what follows unless otherwise stated.

The main emphasis of the variational rationality approach to human behavior is the concept of a {\em worthwhile change}. This approach considers a change to be worthwhile if the motivation to change rather than to stay is higher in comparison with the resistance to change rather than to stay at the current position. To formalize/modelize this concept and related ones, we first clarify the terminology and present the needed definitions.

\subsection{Desirability and Feasibility Aspects of Human Behavior}

Starting with {\bf desirability aspects} of human behavior, we denote by $g\colon X\to\oR:=(-\infty,\infty]$ the so-called ``to be increased" payoffs that reflect profit, satisfaction, utility, valence, revenue functions and the like. ``To be decreased" payoffs are denoted by $\ph(x)=-g(x)$ while reflecting unsatisfied needed, disadvantages, losses, etc. Then:

$\bullet$ {\em Advantage to change} (rather than to stay) payoff from $x$ to $y$ is defined by
\begin{eqnarray}\label{A}
A(y/x):=g(y)-g(x)=\ph(x)-\ph(y)\ge 0.
\end{eqnarray}

$\bullet$ {\em Disadvantage to change} payoff from $x$ to $y$ means accordingly that $A(y/x)\le 0$ for the quantity $A(x/y)$ given in \eqref{A}.

$\bullet$ {\em Motivation to change} from $x$ to $y$ is defined by $M(y/x):=U[A(y/x)]$, where $U\colon A\to\R_+$ is a given nonnegative function of $A$ reflecting the attractiveness of change. In this paper we confine ourselves for simplicity to the linear dependence $U(A)=A$.

Proceeding further with {\bf feasibility aspects} of human behavior, define the linear {\em cost of changing} from $x$ to $y$ by
\begin{eqnarray}\label{C}
C(x,y):=\eta\|y-x\|,\;\mbox{ where the quantity }\;\eta=\eta(x)>0
\end{eqnarray}
denotes the cost of changing per unit of the distance $d(x,y)=\|y-x\|$. Since $\eta(y)\ne\eta(x)$ in many practical situations (e.g., while moving up and down the hill), the cost of changing in \eqref{C} is generally {\em asymmetric}, i.e., $C(y,x)\ne C(x,y)$. We suppose in this paper that $C(x,x)=0$ for all $x\in X$, although it might be that $C(x,x)>0$ in more general settings.

$\bullet$ {\em Inconvenience to change} from $x$ to $y$ (rather than to stay at $x$) is defined by
\begin{eqnarray}\label{I}
I(y/x):=C(x,y)-C(x,x)=C(x,y)\in\R_+.
\end{eqnarray}
$\bullet$ {\em Resistance to change} from $x$ to $y$ (rather than to stay at $x$) is described by $R(y/x):=D[I(y/x)]$ via \eqref{I} and the disutility function $D\colon I\to\R_+$. For simplicity we consider here the linear case where $D(I)=I$ for all $I\ge 0$.

\subsection{Worthwhile Desirable and Feasible Moves}

{\bf Balancing desirability and feasibility aspects.} The variational rationality approach considers different payoffs to balance desirability and feasibility aspects of human dynamics. Among the most important payoffs are the following:

$\bullet$ {\em Proximal payoff} functions are defined by
\begin{eqnarray}\label{pp}
P_\xi(y/x):=g(y)-\xi C(x,y)\;\mbox{ and }\;Q_\xi(y/x):=\ph(y)+\xi C(x,y)
\end{eqnarray}
for ``to be increased" and for ``to be decreased" payoffs $g$ and $\ph$, respectively, where the cost of changing $C(x,y)$ is taken from \eqref{C}. The given {\em weight factor} $\xi=\xi(x)\ge 0$ in \eqref{pp} reflects the status quo at $x$ as well as some other aspects of the current situation.

$\bullet$ {\em Worthwhile to change gain} is defined by
\begin{eqnarray}\label{Aw}
A_\xi(y/x):=A(y/x)-\xi I(y/x)\ge 0
\end{eqnarray}
while representing the difference between advantages \eqref{A} and inconveniences \eqref{I} to change rather than to stay. Based on the definitions, \eqref{Aw} can be represented as the difference of proximal payoffs \eqref{pp} by
\begin{eqnarray}\label{Aw1}
A_\xi(y/x)=\big(g(y)-g(x)\big)-\xi\big(C(x,y)-C(x,x)\big)=P_\xi(y/x)-P_\xi(x/x).
\end{eqnarray}

$\bullet$ {\em Not worthwhile to change loss} is defined similarly to \eqref{Aw} via the {\em loss function} $L(y/x):=-A(y/x)$ by
\begin{eqnarray}\label{loss}
L_\xi(y/x):=L(y/x)+\xi I(y/x)\ge 0,
\end{eqnarray}
i.e., it represents the sum of losses and weighted inconveniences to change rather than to stay. It follows from the definitions that
\begin{eqnarray*}
\begin{array}{ll}
&L_\xi(y/x)=\big(\ph(y)-\ph(x)\big)+\xi\big(C(x,y)-C(x,y)\big)\\
&=\big(\ph(y)+\xi C(x,y)\big)-\big(g(x)+\xi C(x,x)\big)=Q_\xi(y/x)-Q_\xi(x/x),
\end{array}
\end{eqnarray*}
which verifies the relationship $L_\xi(y/x)=-A_\xi(y/x)$ between functions \eqref{Aw} and \eqref{loss}.\\[1ex]
{\bf Worthwhile changes.} Based on the above discussions, we say that it is {\em worthwhile to change} from $x$ to $y$ rather than to stay at the current position $x$ if advantages to change rather than to stay are ``high enough" relatively to inconveniences to change rather than to stay. Mathematically this is expressed via the quantities
\begin{eqnarray}\label{wc}
A(y/x)\ge\xi I(y/x),\;\mbox{ or equivalently as }\;A_\xi(y/x)\ge 0,
\end{eqnarray}
with an appropriate weight factor $\xi=\xi(x)>0$ in \eqref{wc}.\vspace*{0.05in}

Note that the proximal payoff functions \eqref{pp} may be treated as particular cases of the so-called {\em reference-dependent payoffs} $\Gamma(\cdot/r)\colon x\in X\mapsto\Gamma(x/r),\;r\in\R$, which depend on a chosen reference point $r$ even if the original payoffs $g(\cdot)$ and $\ph(\cdot)$ do not depend on $r$.  Observe also that the variational rationality approach addresses also the situations where the original payoff function depend on references points; see \cite{s09,s10} for more details. The latter case covers reference-dependent utilities that were investigated, in particular, in psychology and economics; see, e.g., \cite{kt79,tk91} and the bibliographies therein.

\section{Evaluation Aspects of Variational Rationality}\sce

Due to the extreme complications and numerous uncertain factors of human behavior, reliable {\em evaluation procedures} play a crucial role in behavioral science modeling. The variational rationality approach married to the constructions and results of variational analysis and generalized differentiation offer efficient tools to provide such evaluations.

We start with the relevant evaluation notions from both behavioral and mathematical viewpoints. Roughly speaking, an agent is said to be {\em bounded  rational} if he/she merely tries to improve (``improve enough") his/her ``to be increased" payoff or tries to decrease (``decrease enough") his/her ``to be decreased" payoff. Optimization refers to {\em perfect rationality}. A bounded rational agent has limited capabilities to evaluate advantages and disadvantages to change if
the agent knows his/her payoff $g(x)$ or $\ph(x)=-g(x)$ only at the current position while does not know the values $g(y)$ or  $\ph(y)=-g(y)$ at $y\ne x$.

\subsection{Optimistic and Pessimistic Evaluations}

Consider a realistic case of human behavior where the payoff function $g(\cdot)$ is known at the current position $x$, but unknown (at least precisely) elsewhere for $y\ne x$. First we introduce the notion of {\em optimistic} evaluations and their {\em linear} version.

\begin{Definition}{\bf(optimistic evaluations for ``to be increased" payoffs).}\label{oe} Given a ``to be increased" payoff function $g\colon X\to\oR$ and a current position $x\in X$, we say that:

{\bf(i)} The function $l(\cdot/x)\colon X\to\oR$ provides a {\sc local optimistic evaluation} of $g(\cdot)$ around $x$ if $l(x/x)=g(x)$ and there is a neighborhood $V(x)$ of $x$ such that
\begin{eqnarray*}
g(y)\le l(y/x)\;\mbox{ for all }\;y\in V(x).
\end{eqnarray*}
By definition \eqref{A} of the advantage to change, the above relationships can be unified as
\begin{eqnarray}\label{E0}
A(y/x)=g(y)-g(x)\le l(y/x)-l(x/x)=:E(y/x)\;\mbox{ for all }\;y\in V(x),
\end{eqnarray}
i.e., $A(y/x)\le E(y/x)$ as $y\in V(x)$ for the evaluation quantity $E(y/x)$ defined in \eqref{E0}.

{\bf(ii)} The optimistic evaluation in {\rm(i)} is {\sc linear} if there is $x^*\in X^*$ such that
\begin{eqnarray}\label{lE}
E(y/x)=E_{x^*}(y/x):=\la x^*,y-x\ra\;\mbox{ for all }\;y\in V(x),
\end{eqnarray}
where $x^*$ provides a {\sc rate of change} in the linear optimistic evaluation.

{\bf(iii)} The evaluations in {\rm(i)} and {\rm(ii)} are {\sc global} if $V(x)=X$ therein.

{\bf (iv)} Optimistic evaluations for to ``be decreased" payoffs are defined similarly with replacing $A(y/x)$ by $L(y/x)=-A(y/x)$ and using $L(y/x)\ge E(y/x)$ instead of \eqref{E0}. 
\end{Definition}

To illustrate the notions in Definition~\ref{oe}, we consider the following example that addresses a realistic situation in behavioral economics.

\begin{Example}{\bf(linear optimistic evaluation in economics).}\label{ex1} {\rm Let $g\colon\R_+\to\R$ be a utility function, which depends on consuming a quantity $x\ge 0$ of goods and is defined by
\begin{eqnarray}\label{util}
g(x):=x-\frac{1}{2}x^2\;\mbox{ for all }\;x\in\R_+.
\end{eqnarray}
We can see that the utility function \eqref{util} attends its maximum at $x=1$ while increasing on the interval $[0,1]$ and decreasing on $[0,\infty)$.

$\bullet$ Take $x=1/2$ as our reference point and suppose that the agent knows the value of the utility function \eqref{util} and its derivative therein:
$$
g(1/2)=1/2\;\mbox{ and }\;\nabla g(1/2)=1-1/2=1/2.
$$
Then the local linear evaluation $E_{x^*}(y/x)$ from \eqref{lE} of the advantage to change $A(y/x)$ around the reference position $x=1/2$ is calculated by
\begin{eqnarray}\label{E1}
E_{x^*}(y/x)=x^*(y-x)=\nabla g(x)(y-x)=(1/2)(y-x).
\end{eqnarray}
If in this case we have $x\le y\le 1$, then the computation shows that $0\le A(y/x)\le E_{x^*}(y/x)$. This gives us an optimistic linear evaluation of the realized gain before moving.

If $x\le y\le 1$ with $x=1/2$ as before, then we have the realized loss $-A:=-A(y/x)=g(x)-g(y)\ge 0$ and compute the expected local loss $-E:=-E(y/x)=(1/2)(x-y)\ge 0$ from \eqref{E1}. Since $-A\ge-E$, the linear evaluation is optimistic in this case as well.

$\bullet$ Consider now another initial position $x=3/2$ of the agent with the same utility function $g$ in \eqref{util}. In this setting we have
$$
g(3/2)=3/8\;\mbox{ and }\;\nabla g(x)=-1/2.
$$
Proceeding as above, we examine the following two cases:

Let $1\le y\le x$. Then $A=g(y)-g(x)\ge 0$ is the realized gain and $E=(-1/2)(y-x)\ge 0$ is the expected gain. Since $A\le E$, we have a linear local optimistic evaluation.

Finally, let $x=3/2\le y$. Then we have the realized loss $-A=g(x)-g(y)\ge 0$ and the expected loss $-E=-(-1/2)(y-x)\ge 0$ in this case. Since $-A\ge-E$, it gives us again a linear local optimistic evaluation.

Note that similar justifications can be done for the inequality $L(y/x)\ge E(y/x)$, which describes linear optimistic evaluations of ``to be decreased" payoffs.

Summarizing, we see that in all the cases under considerations the realized gains are {\em lower} than the expected gains while the realized losses are {\em higher} than the expected ones. This is the characteristic feature of {\em optimistic} evaluations.}
\end{Example}

Next we introduce the notion of {\em pessimistic} evaluations where, in contrast to optimistic ones, realized gains are {\em higher} than expected gains while realized losses are {\em lower} than expected ones around the reference point.

\begin{Definition}{\bf (pessimistic evaluations).}\label{pe} Given a payoff function $g\colon X\to\oR$ and a current position $x\in X$, we say that the function $l(\cdot/x)\colon X\to\oR$ provides a {\sc local pessimistic evaluation} of $g(\cdot)$ around $x$ if the following two conditions hold: $l(x/x)=g(x)$ and there is a neighborhood $V(x)$ of $x$ such that
\begin{eqnarray*}
g(y)\ge l(y/x)\;\mbox{ for all }\;y\in V(x).
\end{eqnarray*}
This results in $A(y/x)\ge E(y/x)$ for such $y$. The notions of linear and global evaluations are formulated similarly to cases {\rm(ii)} and {\rm(iii)} of Definition~{\rm\ref{oe}}.
\end{Definition}

Let us now discuss some connections of the above evaluation concepts and their variational interpretations via relevant aspects of optimization and subdifferential theory.\vspace*{0.05in}

$\bullet$ {\bf MM procedures in optimization.} Procedures of this kind and the corresponding algorithms have been well recognized in optimization theory and applications; see, e.g., \cite{bp16} and the references therein. Comparing the general scheme of such algorithms with the above Definitions~\ref{oe} and \ref{pe} of optimistic and pessimistic evaluations shows that {\em majorization-minimization} procedures to minimize a cost function $g\colon X\to\oR$ agree with {\em optimistic} evaluation of $g(y)-g(x)$, while {\em minorization-maximization} procedures to maximize $g(\cdot)$ agree with {\em pessimistic} evaluations of the difference $g(y)-g(x)$.

These observations and the corresponding results of MM optimization allow us to deduce, in particular, that optimistic (resp.\ pessimistic) evaluations of convex (resp.\ concave) functions can always be chosen to be {\em global} and {\em quadratic}.\vspace*{0.05in}

$\bullet$ {\bf Subgradient evaluations of payoffs.} We start discussing here the usage of subgradients of convex and variational analysis to provide constructive linear evaluations of payoff functions. Let us consider for definiteness optimistic evaluations while observing that we can proceed symmetrically with pessimistic ones for symmetric classes of functions (e.g., for concave vs.\ convex). Following the notation and terminology of variational analysis and optimization \cite{m06,rw}, the reference point is denoted by $\ox$ while the moving one by $x$. Given $\vt\colon X\to\oR$ finite at $\ox$ (i.e., with $\ox\in\dom\vt$ from the domain of $\ox$), a {\em classical subgradient} $x^*\in\partial\vt(\ox)$ of $\vt$ at $\ox$ is defined by
\begin{eqnarray}\label{csub}
\la x^*,x-\ox\ra\le\vt(x)-\vt(\ox)\;\mbox{ for all }\;x\in V(\ox),
\end{eqnarray}
where $V(\ox)$ is some neighborhood of $\ox$, and where $\partial\vt(\ox)$ stands for the (local) {\em subdifferential} of $\vt$ at $\ox$ as the collection of all its subgradients at this point. If $V(\ox)=X$ in \eqref{csub}, we get back to the classical definition of the subdifferential of convex analysis, which is nonempty under mild assumptions on a convex function $\vt$; e.g., when $\ox$ is an interior point of $\dom\vt$ (or belongs to the relative interior of $\dom\vt$ in finite dimensions).

\begin{Proposition}{\bf (subgradient linear optimistic evaluations of ``to be decreased" payoffs).}\label{lin-conv} Let $\ph\colon X\to\oR$ be a ``to be decreased" payoff function, let $\ox\in\dom\ph$ be a reference current position, and let
\begin{eqnarray}\label{loss}
L(x/\ox):=\ph(x)-\ph(\ox)\;\mbox{ with }\;x\in V(\ox)
\end{eqnarray}
be the loss function when moving from $\ox$ to the new position $x$ in a neighborhood $V(\ox)$ of $\ox$. Then $x^*\in X^*$ provides a rate of change in the local linear optimistic evaluation of the payoff function $\ph$ around $\ox$ if and only if $x^*\in\partial\ph(\ox)$ in \eqref{csub}. Furthermore, this gives us a complete subdifferential description of the global linear optimistic evaluation of $\ph$ around any $\ox\in\dom\ph$ when $V(\ox)=X$ in \eqref{csub} as in the case of convex functions $\ph(\cdot)$.
\end{Proposition}
{\bf Proof.} If follows from the definitions of subgradients \eqref{csub} and loss function \eqref{loss} that
\begin{eqnarray}\label{ce}
L(x/\ox)\ge\la x^*,x-\ox\ra\;\mbox{ for all }\;x\in V(\ox)
\end{eqnarray}
with the same neighborhood $V(\ox)$ in \eqref{csub} and \eqref{loss}. Invoking now definition \eqref{lE} of the linear quantity $E(x/\ox)$ tells us that
\begin{eqnarray*}
L(x/\ox)\ge E(x/\ox)\;\mbox{ for all }\;x\in V(\ox),
\end{eqnarray*}
which means by Definition~\ref{oe}(iv) that we get a local linear optimistic evaluation of the ``to be decreased" payoff $\ph(\cdot)$ with the rate of change $x^*$. It clearly follows from the above that this evaluation is global if $V(\ox)=X$. $\h$\vspace*{0.05in}

In the next subsection we obtain linear optimistic evaluations of {\em proximal payoffs} associated with general {\em nonconvex} functions via appropriate subgradient extensions.

\subsection{Optimistic Evaluations of Reference-Dependent Payoffs}

Here we address evaluation problems for {\em reference-dependent payoffs} $\Gamma(\cdot/\ox)\colon x\in X\mapsto\Gamma(x/\ox)\in\oR$ discussed in Subsection~2.2. They include, in particular, the proximal payoff functions \eqref{pp} of our special interest in this paper. Similarly to the setting of Definition~\ref{oe}, we say that $x^*\in X^*$ is the {\em rate of change} of the local linear optimistic evaluation for the reference-dependent loss function $L_\Gamma(x/\ox):=\Gamma(x/\ox)-\Gamma(\ox,\ox)$ if there is a neighborhood $V(\ox)$ of the reference point $\ox$ such that
\begin{eqnarray}\label{re}
L_\Gamma(x/\ox)\ge E_{x^*}(x/\ox):=\la x^*,x-\ox\ra\;\mbox{ for all }\;x\in V(\ox).
\end{eqnarray}
It is clear that $x^*$ in \eqref{re} can be interpreted as a (local) classical/convex subgradient \eqref{csub} of the reference-dependent payoff function $\Gamma(\cdot/\ox)$ at $\ox$. Consider now the special case of the reference-dependent {\em proximal payoff function}
\begin{eqnarray}\label{pp-loss}
\Gamma(x/\ox)=Q_\xi(x/\ox):=\ph(x)+\xi\|x-\ox\|,\quad\xi\ge 0,
\end{eqnarray}
defined for the ``to be decreased" original payoff $\ph(\cdot)$ in accordance with \eqref{pp} by taking into account the expression for the changing cost $C(\ox,x)$ in \eqref{C}. Then the local linear optimistic evaluation \eqref{re} is written as
\begin{eqnarray}\label{re1}
L_\xi(x/\ox):=Q_\xi(x/\ox)-Q_\xi(\ox/\ox)\ge E_{x^*}(x/\ox):=\la x^*,x-\ox\ra\;\mbox{ for all }\;x\in V(\ox),
\end{eqnarray}
where $L_\xi(x/\ox)$ can be interpreted as a ``not worthwhile to change loss function." We show next, based on \eqref{re1} and the structure of the proximal payoff $Q_\xi(x/\ox)$ in \eqref{pp-loss}, that the rate of change $x^*$ in the linear optimistic evaluation \eqref{re} can be completely characterized via $\ve$-subgradients of the {\em arbitrary original payoff} $\ph(\cdot)$ in contrast to the convex-like subgradients of the proximal payoff function $Q_\xi(x/\ox)$ as in \eqref{re1}.

Given an arbitrary function $\vt\colon X\to\oR$ finite at $\ox$ and following \cite[Definition~1.83]{m06}, we say that $x^*$ is an {\em $\ve$-subgradient} of $\vt$ at $\ox$ for some $\ve\ge 0$ if it belongs to the (analytic) {\em $\ve$-subdifferential} of $\vt$ at this point defined by
\begin{eqnarray}\label{e-sub}
\Hat\partial_\ve\vt(\ox):=\Big\{x^*\in X^*\Big|\;\disp\liminf_{x\to\ox}\frac{\vt(x)-\vt(\ox)-\la x^*,x-\ox\ra}{\|x-\ox\|}\ge-\ve\Big\}.
\end{eqnarray}

The next theorem establishes a {\em two-sided relationship} between local linear optimistic evaluations of ``to be decreased" payoff functions in the variational rationality theory for behavioral sciences and $\ve$-subgradients of variational analysis. From one side it reveals a {\em behavioral sense} of $\ve$-subgradients for general functions, while from the other side it provides an efficient $\ve$-{\em subdifferential mechanism} to calculate rates of change in local linear optimistic evaluations of ``to be decreased" payoffs. Note that the second assertion of the theorem allows us to keep a similar description with the replacement of the original payoff by its differentiable approximation with the same rate of change at the reference point. This makes the result to be easier for implementations and applications.

\begin{Theorem}{\bf ($\ve$-subgradient description of rates of change in linear optimistic evaluations of ``to be decreased" payoffs).}\label{e-desc} Let $\ph\colon X\to\oR$ be a ``to be decreased" payoff in behavioral dynamics, and let $\ve\ge 0$.

{\bf (i)} We have that $x^*\in\Hat\partial_\ve\ph(\ox)$ if and only if for every weight factor $\xi>\ve$ in the proximal payoff \eqref{pp-loss} there is a neighborhood $V(\ox)$ of $\ox$ such that the local linear optimistic evaluation \eqref{re1} holds around $\ox$ with the rate of change $x^*$ for the not worthwhile to change loss function $L_\xi(x/\ox)$ represented as
\begin{eqnarray}\label{sm}
L_\xi(x/\ox)=\big(\ph(x)-\ph(\ox)\big)+\xi\|x-\ox\|.
\end{eqnarray}

{\bf (ii)} It $\ve=0$, then in addition to the behavioral description of $x^*\in\Hat\partial_0\ph(\ox):=\Hat\partial\ph(\ox)$ in {\rm(i)} we claim the existence of a ``more decreased" function $s\colon V(\ox)\to\R$, which is Fr\'echet differentiable at $\ox$ with $\nabla s(\ox)=x^*$, $s(\ox)=\ph(\ox)$, and $s(x)\le\ph(x)$ whenever $x\in V(\ox)$.
\end{Theorem}
{\bf Proof.} To justify (i), we first fix $\ve\ge 0$ and pick an arbitrary $\ve$-subgradient $x^*\in\Hat\partial_\ve\ph(\ox)$ for the ``to be decreased" payoff $\ph(\cdot)$ under consideration. Given $\xi>\ve$, denote $\nu:=\xi-\ve>0$ and employ the {\em variational description} of $\ve$-subgradients from \cite[Proposition~1.84]{m06}, which says that for every $\nu>0$ there is a neighborhood $V(\ox)$ of $\ox$ such that the function
\begin{eqnarray*}
\psi(x):=\ph(x)-\ph(\ox)+(\ve+\nu)\|x-\ox\|-\la x^*,x-\ox\ra
\end{eqnarray*}
attains it minimum on $V(\ox)$ at $\ox$. Thus $\psi(x)\ge\psi(\ox)=0$ whenever $x\in V(\ox)$. This yields
\begin{eqnarray}\label{sm1}
\ph(x)-\ph(\ox)+\xi\|x-\ox\|\ge\la x^*,x-\ox\ra\;\mbox{ for all }\;x\in V(\ox)
\end{eqnarray}
by taking into account the choice of $\nu$. Recalling the definition of the proximal payoff $Q_\xi(x/\ox)$ in \eqref{pp-loss} with $Q_\xi(\ox/\ox)=\ph(\ox)$, we deduce from \eqref{sm1} and \eqref{re1} that $x^*$ is a rate of change in the local linear optimistic evaluation of $L_\xi(x/\ox)$ from \eqref{sm} around $\ox$.

To verify the converse implication in (i), we simply reverse the argumentations above while observing that the variational description of $\ve$-subgradients in \cite[Proposition~8.4]{m06} provides a complete characterization of these elements.

It remains to justify assertion (ii) of the theorem. To proceed, we employ a {\em smooth variational description} of {\em regular subgradients} $x^*\in\Hat\partial\ph(\ox)$ (known also as Fr\'echet and viscosity ones) from \cite[Theorem~1.88]{m06}(i), which ensures the existence of a real-valued function $s\colon V(\ox)\to\R$ satisfying the listed properties. The aforementioned result is stated in Banach spaces $X$ while its proof holds without the completeness requirement on the normed space $X$ under consideration. $\h$

\begin{Remark}{\bf (further results and discussions on linear subgradient evaluations).}\label{rem1} {\rm

{\bf (i)} The properties of the supporting function $s(\cdot)$ in Theorem~\ref{e-desc} can be significantly improved if the space $X$ is Banach and satisfies additional ``smoothness" requirements; see \cite[Theorem~1.88(ii,iii)]{m06} for more details.

{\bf (ii)} It has been well recognized in variational analysis that $\ve$-subgradients ($\ve\ge 0$) from \eqref{e-sub} used in linear optimistic evaluations of Theorem~\ref{e-desc} are not robust and do not satisfy major calculus rules for subdifferentiation of sums, compositions, etc., unless $\vt$ is convex and belong to some other restrictive classes of functions. This may complicate applications of subgradient linear evaluations established in Theorem~\ref{e-desc}. The situation is dramatically improved for the limiting construction
\begin{eqnarray}\label{ls}
\begin{array}{ll}
\partial\vt(\ox):=\Big\{x^*\in X^*\Big|&\exists\,\mbox{ seqs. }\;\ve_k\dn 0,\;x_k\st{\vt}{\to}\ox,\;\mbox{ and }\;\;x^*_k\st{w^*}{\to}x^*\\
&\mbox{such that }\;x^*_k\in\Hat\partial_{\ve_k}\vt(x_k)\;\mbox{ for all }\;k=1,2,\ldots\Big\}
\end{array}
\end{eqnarray}
known as the Mordukhovich {\em basic/limiting subdifferential} of $\ph$ at $\ox\in\dom\vt$. The symbol $x\st{\vt}{\to}\ox$ in \eqref{ls} indicates that $x\to\ox$ with $\vt(x)\to\vt(\ox)$ while $w^*$ stands for the weak$^*$ topology of $X^*$. The limiting subdifferential construction \eqref{ls} can be viewed as a {\em robust regularization} of $\ve$-subdifferentials \eqref{e-sub}, and the validity of {\em full calculus} for it in finite-dimensions and broad infinite-dimensional settings (particularly in the class of {\em Asplund} spaces, which includes every reflexive Banach space, etc.) is due to {\em variational/extremal principles} of variational analysis; see \cite{m06} for a comprehensive study. Having this in mind, we may treat robust subdifferential evaluations of behavioral payoffs expressed in terms of \eqref{ls} as {\em asymptotic} versions of those obtained in Theorem~\ref{e-desc} and in what follows.}
\end{Remark}

\section{Exact Stationary Traps in Behavioral Dynamics}\sce

One of the most important questions in the framework of variational rationality can be formulated as follows: {\em when is it worthwhile to change the current position rather than to stay at it?} This issue is closely related to appropriate notions of {\em traps}. Among the main aims of this paper is to study some versions of {\em stationary traps} and their efficient descriptions by using adequate tools of variational analysis and generalized differentiation. We start our study by considering {\em exact} versions of stationary traps in behavioral models and then proceed with their {\em approximate} counterparts in the next section. To begin with, we formally introduce the notions under consideration.

\begin{Definition}{\bf (stationary traps).}\label{traps} Let $\ox\in X$ be a reference point, let $\xi\ge 0$ be a weight factor, and let $Q_\xi(x/\ox):=\ph(x)+\xi\|x-\ox\|$ be the reference-dependent proximal payoff built upon a given ``to be decreased" original payoff $\ph\colon X\to\oR$. We say that:

{\bf (i)} $\ox$ is a {\sc local stationary trap} in behavior dynamics with the weight factor $\xi$ if there is a neighborhood $V(\ox)$ such that
\begin{eqnarray}\label{loc-trap}
L_\xi(x/\ox):=Q_\xi(x/\ox)-Q_\xi(\ox/\ox)\ge 0\mbox{ for all }\;x\in V(\ox)
\end{eqnarray}
via the not worthwhile to change loss function $L_\xi(x/\ox)$ that can be directly defined by \eqref{sm}. Equivalently, condition \eqref{loc-trap} can be expressed in the form
\begin{eqnarray}\label{trap1}
A_\xi(x/\ox):=A(x/\ox)-\xi I(x/\ox)\le 0\;\mbox{ for all }\;x\in V(\ox),
\end{eqnarray}
where $A(x/\ox)$ and $I(x/\ox)$ are the advantage to change and the inconvenience to change from $\ox$ to $x$ defined in terms of the ``to be increased" payoff function $g(x)=-\ph(x)$ in \eqref{A} and \eqref{I}, respectively, via the cost of changing \eqref{C} with $\eta=1$ for simplicity.

{\bf (ii)} $\ox$ is a {\sc strict local stationary trap} in behavioral dynamics with the weight factor $\xi$ if there exists a neighborhood $V(\ox)$ such that a counterpart of \eqref{loc-trap} holds with the replacement of ``$\ge$" by ``$>$" for all $x\in V(\ox)\setminus\{\ox\}$.

{\bf (iii)} If $V(\ox)=X$, then the conditions in {\rm(i)} and {\rm(ii)} define $\ox$ as a {\sc global stationary trap} and its {\sc strict global} version, respectively.
\end{Definition}

Roughly speaking, stationary traps are positions, which are {\em not worthwhile to quit}. For definiteness, we confine ourselves in what follows to considering only stationary traps while observing that their strict counterparts can be studied similarly.\vspace*{0.05in}

The subdifferential optimistic evaluations of reference-dependent proximal payoffs obtained in Theorem~\ref{e-desc} allow us to derive efficient {\em certificates} (sufficient conditions) for stationary traps in the general framework of normed spaces $X$. Again for definiteness, the results below are presented only in terms of ``to be decreased" payoff functions.

\begin{Proposition}{\bf ($\ve$-subdifferential certificates for stationary traps from optimistic evaluations).}\label{cert-traps} Given a ``to be decreased" payoff $\ph\colon X\to\oR$ and a point $\ox\in\dom\ph$, fix any $\ve\ge 0$ and assume that there exist an $\ve$-subgradient $x^*\in\Hat\partial_\ve\ph(\ox)$ and a neighborhood $V(\ox)$ of $\ox$ such that we have
\begin{eqnarray}\label{str-trap}
E_{x^*}(x/\ox):=\la x^*,x-\ox\ra\ge 0\;\mbox{ for all }\;x\in V(\ox).
\end{eqnarray}
Then $\ox$ is a local stationary trap in behavioral dynamics with any weight factor $\xi>\ve$. Furthermore, this trap is global if $V(\ox)=X$ in \eqref{str-trap}.
\end{Proposition}
{\bf Proof.} It follows from Theorem~\ref{e-desc}(i) that for any $\ve\ge 0$, any $\ve$-subgradient $x^*\in\Hat\partial_\ve\ph(\ox)$, and any weight $\xi>\ve$ there is a neighborhood $V(\ox)$ such that
\begin{eqnarray}\label{e-trap}
L_\xi(x/\ox):=Q_\xi(x/\ox)-Q_\xi(\ox/\ox)\ge\la x^*,x-\ox\ra\;\mbox{ for any }\;x\in V(\ox)
\end{eqnarray}
Combining \eqref{e-trap} with assumption \eqref{str-trap} and using the stationary trap Definition~\ref{traps}(i,iii) justifies the claimed conclusions. $\h$\vspace*{0.05in}

Since we clearly have by definition \eqref{e-sub} that
\begin{eqnarray*}
\Hat\partial_{\ve_1}\ph(\ox)\subset\Hat\partial_{\ve_2}\ph(\ox)\;\mbox{ whenever }\;0\le\ve_1\le\ve_2,
\end{eqnarray*}
for each $\ox$ and $x^*$ there is the {\em minimal subdifferential factor} $\ve_{\rm\small min}=\ve_{\rm\small min}(\ox,x^*)=\min\{\ve\}$ over all $\ve$ such that \eqref{str-trap} holds. Note that we may have $\ve_{\rm\small min}>0$ for some $\ox$ and $x^*$.

It follows from the subgradient estimate \eqref{str-trap} that the analysis of stationary traps via Proposition~\ref{cert-traps} depends on the two major parameters:\vspace*{0.02in}

{\bf (a)} an $\ve$-subgradient $x^*$ that can be chosen from the subdifferential set $\Hat\partial_\ve\ph(\ox)$ together with an appropriate number $\ve\ge 0$;

{\bf (b)} a weight factor $\xi>\ve$ that is generally unknown, or at least not known exactly. \vspace*{0.05in}

Let us discuss some situations that emerge from Proposition~\ref{cert-traps} when $\xi>\ve$ and also directly from Definition~\ref{traps} when the latter condition fails.\vspace*{0.05in}

{\bf 1. Flat linear optimistic evaluations}. This is the case when there exists $\ve\ge 0$ such $0\in\Hat\partial_\ve\ph(\ox)$. Choosing $x^*=0$, we get the flat evaluation
\begin{eqnarray*}
E_{x^*}(x/\ox)=\la x^*,x-\ox\ra=0\;\mbox{ for all }\;x\in V(\ox),
\end{eqnarray*}
which ensures that $\ox$ is a local/global stationary trap by Proposition~\ref{cert-traps}.

$\bullet$ If $\xi=0$ in this setting, Proposition~\ref{cert-traps} does not apply although $E_{x^*}(x/\ox)=0$. However, the very definitions of local (global) traps in \eqref{loc-trap} amount to saying that
\begin{eqnarray*}
\ph(x)-\ph(\ox)\ge 0\;\mbox{ for all }\;x\in V(\ox),
\end{eqnarray*}
which means that $\ox$ is a local (global) {\em minimizer} of the {\em original payoff} $\ph(\cdot)$. Thus we have in this case that local (global) stationary traps in  behavioral dynamics correspond to local (global) minima of ``to be decreased" payoff functions.

$\bullet$ If $\xi>0$ in this setting, then stationary traps mean by \eqref{loc-trap} that
\begin{eqnarray*}
Q_\xi(x/\ox)-Q_\xi(\ox/\ox)\ge 0\;\mbox{ for all }\;x\in V(\ox),
\end{eqnarray*}
i.e., local (global) stationary traps reduce in this to local (global) {\em minimizers} of the {\em proximal payoff} $Q(x/\ox)$ instead of the original one independently of $\xi>0$. If $\xi>\ve$, it can be also detected via the linear optimistic evaluation of Proposition~\ref{cert-traps}.\vspace*{0.05in}

{\bf 2. Subdifferential evaluations of variational analysis.} This is the setting where $x^*\ne 0$ for a given $\ve$-subgradient from the set $\Hat\partial_\ve\ph(\ox)$ as $\ve\ge 0$.

$\bullet$ If $\xi=0$ (no inertia and costs of changing in the moving process), Proposition~\ref{cert-traps} does not apply, while the application of \eqref{loc-trap} is exactly the same as in the case of $x^*=0$.

$\bullet$ If $\xi>0$, then there are inertia and costs of changing in the moving process, while the application of \eqref{loc-trap} is not different from the case of $x^*=0$ and does not actually provide verifiable information for the stationary trap determination. This is due to the fact that the very definition \eqref{loc-trap} with $\xi>0$ is merely conceptional, not constructive, and that the weight factor $\xi$ is generally {\em unknown}. In contrast to it, the subgradient linear optimistic evaluation \eqref{str-trap} obtained in Proposition~\ref{cert-traps} ensures that the position $\ox$ is a local or global trap {\em for any weight factor} $\xi>\ve$ without determining this factor a priori. This is a {\em serious advantage} of the subgradient trap evaluation from variational analysis.

$\bullet$ The relationship $\xi>\ve$  between the weight factor $\xi$ (a parameter of {\em behavioral dynamics}) and the subdifferential factor $\ve$ (a parameter of {\em variational analysis}) plays a crucial role in the obtained evaluation of stationary traps. Let us present a striking interpretation of this relationship from the viewpoint of {\em variational stationarity}.

It is reasonable to identify the subdifferential factor $\ve$ with the {\em unit cost of changing} $\eta=\eta(\ox)$ in $C(\ox,x)$ from $\ox$ to $x$ in \eqref{C}. This means that the value $\ve=\eta$ reflects the {\em resistance aspect} of behavioral dynamics near $\ox$. The lower number $\ve$, the less resistance of the agent to change locally around the reference point is. On the other hand, the size of $\xi$ determines how much it is {\em worthwhile to move} from $\ox$ to $x$. Thus the required condition $\xi>\ve$ shows that the advantage to change should be {\em high enough} in comparison to the resistance to change at the reference position of the agent.

\section{Evaluations of Approximate Stationary Traps}\sce

In this subsection we consider more flexible {\em approximate} versions of stationary traps from Definition~\ref{traps} (omitting their strict counterparts) and show that these notions admit verifiable subdifferential evaluations similar to those obtained above for exact traps as well as new variational descriptions derived by using fundamental {\em variational principles}.

\begin{Definition}{\bf (approximate stationary traps).}\label{atraps} Staying in the framework of Definition~{\rm\ref{traps}}, take any $\gg>0$. It is said that $\ox$ is a {\sc $\gg$-approximate local stationary trap} in behavioral dynamics with the weight factor $\xi\ge 0$ if there is a neighborhood $V(\ox)$ of $\ox$ on which we have the inequality
\begin{eqnarray}\label{app-trap}
L_\xi(x/\ox):=Q_\xi(x/\ox)-Q_\xi(\ox/\ox)\ge-\gg\;\mbox{ whenever }\;x\in V(\ox).
\end{eqnarray}
The $\gg$-approximate stationary trap $\ox$ is {\sc global} if $V(\ox)=X$ in \eqref{app-trap}.
\end{Definition}

Similarly to the case of stationary traps we derive efficient certificates of approximate stationary traps from optimistic subdifferential evaluations of
proximal payoffs.

\begin{Proposition}{\bf ($\ve$-subdifferential certificates for approximate stationary traps from optimistic evaluations).}\label{app-cert}  Given a ``to be decreased" payoff $\ph\colon X\to\oR$, a point $\ox\in\dom\ph$ and a rate $\gg>0$, fix any $\ve\ge 0$ and assume that there exist an $\ve$-subgradient $x^*\in\Hat\partial_\ve\ph(\ox)$ and a neighborhood $V(\ox)$ of $\ox$ such that we have
\begin{eqnarray}\label{app-trap1}
E_{x^*}(x/\ox):=\la x^*,x-\ox\ra\ge-\gg\;\mbox{ for all }\;x\in V(\ox).
\end{eqnarray}
Then $\ox$ is a $\gg$-approximate local stationary trap with any weight factor $\xi>\ve$. Furthermore, the $\gg$-approximate stationary trap is global if $V(\ox)=X$ in \eqref{app-trap1}.
\end{Proposition}
{\bf Proof.} As in the proof of Proposition~\ref{cert-traps}, we deduce from Theorem~\ref{e-desc}(i) that estimate \eqref{e-trap} holds for any $\ve\ge 0$, any $\ve$-subgradient $x^*\in\Hat\partial_\ve\ph(\ox)$, and any weight $\xi>\ve$ with some neighborhood $V(\ox)$, which can be identified with the one in \eqref{app-trap1} without loss of generality. Then substituting the assumed estimate \eqref{app-trap1} into \eqref{e-trap} shows that $\ox$ is a $\gg$-approximate stationary trap (local or global) by definition \eqref{app-trap}. $\h$\vspace*{0.05in}

Based on the certificate obtained in Proposition~\ref{app-cert} when $\xi>\ve$ and on definition \eqref{app-trap1} otherwise, we discuss some remarkable features of approximate stationary traps.

Observe first that {\em flat evaluations} of approximate stationary traps is about the same as for their exact counterparts discussed above. Indeed, in this case
we simply replace exact minimizers of the original payoff $\ph(\cdot)$ by its {\em $\gg$-approximate minimizers}
\begin{eqnarray*}
\ph(\ox)\le\ph(x)+\gg\;\mbox{ for all }\;x\in V(\ox)
\end{eqnarray*}
if $\xi=0$, and correspondingly for the proximal payoff $Q_\xi(x/\ox)$ if $\xi>0$.

If $\xi>\ve$, we can apply the {\em subdifferential evaluations} from Proposition~\ref{app-cert} with the estimate $E_{x^*}(x/\ox)$ from \eqref{app-trap1}.
we always have
\begin{eqnarray}\label{E}
\big|E_{x^*}(x/\ox)\big|=\big|\la x^*,x-\ox\ra\big|\le\|x^*\|\cdot\|x-\ox\|,
\end{eqnarray}
which allows us to single out the case of {\em flat enough evaluations} for approximate stationary traps. It is the case when the value of $\|x^*\|$ is {\em sufficiently small} for some $\ve$-subgradient $x^*\in\Hat\partial_\ve\ph(\ox)$. In this case the estimating value $|E_{x^*}(x/\ox)|$ can be made small enough even for a large set $V(\ox)$ in \eqref{app-trap1}. This verifies by Proposition~\ref{app-cert} that there exists $\gg>0$ such that $\ox$ is a $\gg$-approximate stationary trap \eqref{app-trap}

Note that, since the expected gain or loss $E_{x^*}(x/\ox)$ is small when $\|x^*\|$ is small, the agent has no incentive to change from $\ox$ to $x$ even in the case of gain. Changing a bit (a marginal change) leads to {\em disappointment} in both cases of loss and gain. Thus an approximate local stationary trap is a position where it is {\em not worthwhile to change}.\vspace*{0.05in}

Next we proceed with applying a different machinery of variational analysis to study approximate traps by involving powerful {\em variational principles} instead of subdifferential estimates as above. The new machinery essentially relies on the fact that $\gg\ne 0$ in Definition~\ref{atraps} of $\gg$-approximate traps, i.e., it does not apply to the study of exact traps. Furthermore, in contrast to the results of Propositions~\ref{cert-traps} and \ref{app-cert}, the statements below reveal conditions that are satisfied at the approximate trap point $\ox$, i.e., {\em necessary conditions} for this concept. Recall that Proposition~\ref{app-cert} offers {\em sufficient conditions} of a different type that ensure the validity of the approximate trap property \eqref{app-trap}.\vspace*{0.05in}

To proceed in the new direction, we first employ an appropriate version of the fundamental {\em Ekeland variational principle} to establish the {\em existence} of approximate stationary traps with certain {\em optimality} properties to minimize {\em perturbed} proximal payoffs.

\begin{Theorem}{\bf (existence of approximate stationary traps minimizing perturbed proximal payoffs).}\label{var-traps} Pick any $\gg>0$ and let $\ox\in\dom\ph$ be a local $\gg$-approximate stationary trap from Definition~{\rm\ref{atraps}} with the "to be decreased" payoff $\ph\colon X\to\oR$ and the neighborhood $V(\ox)$ in \eqref{app-trap}. Assume that the space $X$ is Banach, and that the function $\ph(\cdot)$ is lower semicontinuous and bounded from below around $\ox$. Then for every $\xi\ge 0$ and every $\lm>0$ there exists $x_\gg=x_\gg(\xi,\lm)$ such that:

{\bf (a)} $x_\gg$ is a local $\gg$-approximate stationary trap;

{\bf (b)} $\|x_\gg-\ox\|\le\lm$;

{\bf (c)} $x_\gg$ is an exact minimizer of the perturbed proximal payoff meaning that
\begin{eqnarray}\label{pert}
Q_\xi(x_\gg/\ox)\le Q_\xi(x/\ox)+\frac{\gg}{\lm}\|x-x_\gg\|\;\mbox{ for all }\;x\in V(\ox).
\end{eqnarray}
\end{Theorem}
{\bf Proof.} Assume without loss of generality that $V(\ox)$ is closed and define the function
\begin{eqnarray}\label{th}
\th(x):=L_\xi(x/\ox)\;\mbox{ for all }\;x\in V(\ox),
\end{eqnarray}
which is lower semicontinuous and bounded from below on the complete metric space $V(\ox)$. It follows from \eqref{sm} and the $\gg$-approximate stationary trap definition \eqref{app-trap} that
\begin{eqnarray*}
\th(\ox)=0\;\mbox{ and }\;\th(\ox)\le\inf_{V(\ox)}\th(x)+\gg\;\mbox{ for all }\;x\in V(\ox).
\end{eqnarray*}
Thus we are in a position to apply the Ekeland variational principle (see, e.g., \cite[Theorem~2.26]{m06}) to $\th$ on $V(\ox)$ with the initial data $x_0:=\ox$, $\ve:=\gg$, and $\lm:=\lm$ taken from the formulation of the theorem. This ensures the existence of a perturbed point $x_\gg\in V(\ox)$ such that $\|x_\gg-\ox\|\le\lm$, that the minimality condition
\begin{eqnarray}\label{pert1}
\th(x_\gg)\le\th(x)+\frac{\gg}{\lm}\|x-x_\gg\|\;\mbox{ for all }\;x\in V(\ox)
\end{eqnarray}
is satisfied, and that $\th(x_\gg)\le\th(\ox)$. The latter condition immediately implies that $x_\gg$ is also a local $\gg$-approximate stationary trap in \eqref{app-trap}, while \eqref{pert1} readily yields the minimality condition \eqref{pert}. This therefore completes the proof of the theorem. $\h$\vspace*{0.05in}

Note that the perturbation term in \eqref{pert} can be made {\em regulated} by an arbitrary choice of $\lm>0$ under the fixed rank $\gg>0$ of both approximate stationary traps $\ox$ and $x_\gg$. Note furthermore that the perturbation structure in \eqref{pert} relates to the {\em changing cost} $C(x_\gg,x)$ from \eqref{C} at the minimizing $\gg$-approximate stationary trap $x_\gg$. Overall we can see the following {\em striking behavioral interpretation} of Theorem~\ref{var-traps}: starting with an arbitrary local $\gg$-approximate trap $\ox$ it is possible to reach in one arbitrary small step $\|x_\gg-\ox\|\le\lm$ another local $\gg$-approximate trap at which we precisely minimize the perturbed proximal payoff function $Q_\xi(\cdot/\ox)$. \vspace*{0.05in}

The following variational result establishes the existence of a local $\gg$-approximate stationary trap in the modified framework of Theorem~\ref{var-traps} with a regulated {\em rate of change} at this trap, which can be made sufficiently small if needed. The proof is based on applying the {\em lower subdifferential variational principle} \cite{m06} that is actually equivalent by \cite[Theorem~2.28]{m06} to the (geometric) {\em approximate extremal principle} considered in the next section. Furthermore, the aforementioned result from \cite{m06} says that a large subclass of Banach spaces, known as Asplund spaces, is the most adequate setting for the validity of these fundamental principles of variational analysis.

Recall that a Banach space $X$ is {\em Asplund} if each of its separable subspaces has a separable dual. This class is sufficiently broad while including, in particular, each reflexive spaces, Banach spaces admitting equivalent norms that are Fr\'echet differentiable at nonzero points, those with separable topological duals, etc. Besides the aforementioned lower subdifferential variational principle, the following theorem exploits the subdifferential {\em sum rules} that are valid for regular subgradients \eqref{e-sub} as $\ve=0$ under the Fr\'echet differentiability of one summand and for our basic/limiting subgradients \eqref{ls} in Asplund spaces without any differentiability assumptions. The latter sum rule fails for the regular subdifferential even in simple nonsmooth settings.

\begin{Theorem}{\bf (subgradient rates of change at local approximate stationary traps).}\label{sub-trap} Staying in the framework of Theorem~{\rm\ref{var-traps}}, assume in addition that the space $X$ is Asplund. Then for every $\xi\ge 0$ and every $\lm>0$ there exists $x_\gg=x_\gg(\xi,\lm)$ such that conditions {\rm(a)} and {\rm(b)} therein hold, while condition {\rm(c)} is replaced by the following:
\begin{eqnarray}\label{rate}
\mbox{there exists }\;x^*_{\gg}\in\Hat\partial Q_\xi(x_\gg/\ox)\;\mbox{ with }\;\|x^*_{\gg}\|\le\frac{\gg}{\lm}.
\end{eqnarray}
Furthermore, we have the following inclusion for the rate of change in \eqref{rate}:
\begin{eqnarray}\label{rate1}
x^*\in\partial\ph(x_\gg)+\left\{\begin{array}{ll}
\big\{v^*\in X^*\big|\;\|v^*\|\le 1\big\}\;\mbox{ if }\;x_\gg=\ox,\\\\
\big\{v^*\in X^*\big|\;\|v^*\|=1,\;\la v^*,x_\gg-\ox\ra=\|x_\gg-\ox\|\big\}\;\mbox{ if }\;x_\gg\ne\ox.
\end{array}
\right.
\end{eqnarray}
Finally, the inclusion for the rate $x^*$ reads as
\begin{eqnarray}\label{rate2}
x^*\in\nabla\ph(x_\gg)+\left\{\begin{array}{ll}
\big\{v^*\in X^*\big|\;\|v^*\|\le 1\big\}\;\mbox{ if }\;x_\gg=\ox,\\\\
\big\{v^*\in X^*\big|\;\|v^*\|=1,\;\la v^*,x_\gg-\ox\ra=\|x_\gg-\ox\|\big\}\;\mbox{ if }\;x_\gg\ne\ox
\end{array}
\right.
\end{eqnarray}
provided that the original payoff $\ph(\cdot)$ is Fr\'echet differentiable at $\ox$.
\end{Theorem}
{\bf Proof.} To proceed, we consider the lower semicontinuous and bounded from below function $\th(\cdot)$ from \eqref{th} defined on an Asplund space and apply to it the lower subdifferential variational principle from \cite[Theorem~2.28]{m06} with the parameters $x_0:=\ox$, $\ve:=\gg$, and $\lm:=\lm$. This gives us a perturbed point $x_\gg\in V(\ox)$ satisfying $\|x_\gg-\ox\|\le\lm$ and $\th(x_\gg)\le\th(\ox)$ as well as a regular subgradient $x^*_\gg\in\Hat\partial\th(x_\gg)$ with $\|x^*_\gg\|\le\gg/\lm$. The latter inclusion clearly gives us \eqref{rate} due to \eqref{th} and the form of $L_\xi(x/\ox)$ in \eqref{app-trap}.

To justify further the subdifferential inclusion \eqref{rate} with our basic subdifferential \eqref{ls} therein, we first observe from \eqref{ls} that
\begin{eqnarray*}
\Hat\partial Q_\xi(x_\gg/\ox)\subset\partial Q_\xi(x_\gg/\ox).
\end{eqnarray*}
Then we apply the basic subdifferential sum rule from \cite[Theorem~2.33]{m06} to the proximal payoff function defined as the sum
\begin{eqnarray*}
Q_\xi(x/\ox)=\ph(x)+\xi\|x-\ox\|,\quad\xi\ge 0
\end{eqnarray*}
at $x=x_\gg$ with taking into account that the function $\xi\|x-\ox\|$ is Lipschitz continuous around $\ox$ while the original payoff $\ph(\cdot)$ is assumed to be lower semicontinuous around this point. Using the sum rule and the well-known formula for subdifferentiation of the norm function in convex analysis, we deduce from \eqref{rate} the subdifferential inclusion \eqref{rate1}.

It remains to verify the final rate of change inclusion \eqref{rate2} when $\ph(\cdot)$ is Fr\'echet differentiable at $x_\gg$. It is in fact a direct consequence of the equality sum rule from \cite[Proposition~1.107(i)]{m06} applied to $\Hat\partial Q_\xi(x_\gg/\ox)$ in \eqref{rate} due to the summation form of $Q_\xi(\cdot/\ox)$. This completes the proof of the theorem. $\h$\vspace*{0.05in}

Note that the last inclusion \eqref{rate2} is generally independent of the subdifferential one \eqref{rate1} since $\partial\ph(x_\gg)$ may not reduce to the gradient $\nabla\ph(x_\gg)$ when $\ph(\cdot)$ is merely Fr\'echet differentiable at this point. A simple example is provided by the function $\ph(x):=x^2\sin(1/x)$ if $x\ne 0$ and $\ph(0):=0$ with $x_\gg=0$. For this function we have
\begin{eqnarray*}
\nabla\ph(0)=0\;\mbox{ while }\;\partial\ph(0)=[-1,1].
\end{eqnarray*}
For the validity of $\partial\ph(x_\gg)=\{\nabla\ph(x_\gg)\}$ we need to require the {\em strict differentiablity} of $\ph(\cdot)$ at the point in question; this holds, in particular, when $\ph(\cdot)$ continuously differentiable around this point; see \cite[Definition~1.13 and Corollary~1.82]{m06}.

Observe finally that in the case where $x_\gg\ne\ox$ and the norm on $X$ is Fr\'echet differentiable at nonzero points (we can reduce to this setting any Banach space with a Fr\'echet differentiable renorm; in particular, any reflexive space \cite{m06}), the second summand in \eqref{rate1} and \eqref{rate2} is a {\em singleton}. Thus in this case formula \eqref{rate2} holds as equality, and we get a {\em precise computation} of the rate of change $x^*$ at the $\gg$-approximate stationary trap $x_\gg$.

\section{Geometric Evaluations and Extremal Principle}\sce

In this section we continue the investigation of stationary traps in behavioral dynamics while invoking now for these purposes {\em geometric} constructions and results of variational analysis. Given a set $\O\subset X$, a point $\ox\in\O$, and a number $\ve\ge 0$, we consider the collection of $\ve$-{\em normals} to $\O$ at $\ox$ defined by
\begin{eqnarray}\label{nor}
\Hat N_\ve(\ox;\O):=\Big\{x^*\in X^*\Big|\;\disp\limsup_{x\st{\O}{\to}\ox}\frac{\la x^*,x-\ox\ra}{\|x-\ox\|}\le\ve\Big\},
\end{eqnarray}
where the symbol $x\st{\O}{\to}\ox$ signifies that $x\to\ox$ with $x\in\O$. Let $u_{x^*}(x):=\la x^*,x\ra$ be a linear {\em utility function} with the rate of change $x^*$, and let $\xi\ge 0$. Define the {\em proximal payoff} of type \eqref{pp} for to ``to be increased" utility function $u_{x^*}(x)$ by
\begin{eqnarray}\label{ppg}
P_{\xi,x^*}(x/\ox):=u_{x^*}(x)-\xi\|x-\ox\|.
\end{eqnarray}

The following result gives us a generalized normal description of local stationary traps in behavioral dynamics with linear utility functions. We say that $\ox$ is a local stationary trap {\em relative to} some set $\O$ if $x$ belongs to $\O$ in the relationships of Definition~\ref{traps}(i). Our study below is based on the "to be increased" payoff form \eqref{trap1} of stationary traps.

\begin{Proposition}{\bf ($\ve$-normal description of local stationary traps).}\label{nor-trap} Let $\ox\in\O$, and let $\ve\ge 0$. Then $x^*\in\Hat N_\ve(\ox;\O)$ if and only if $\ox$ is a local stationary trap relative to $\O$ with the linear utility function $u_{x^*}(x)=\la x^*,x\ra$ and any weight factor $\xi>\ve$.
\end{Proposition}
{\bf Proof.} We employ the {\em variational description} of $\ve$-normals from \cite[Proposition~1.28]{m06} saying that $x^*\in\Hat N_\ve(\ox;\O)$ if and only if for any $\nu>0$ the function
\begin{eqnarray*}
\psi(x):=\la x^*,x-\ox\ra-(\ve+\nu)\|x-\ox\|
\end{eqnarray*}
attains its local maximum $\psi(\ox)=0$ at $\ox$ over $\O$. This means that there is a neighborhood $V(\ox)$ of $\ox$ relative to $\O$ such that the proximal payoff \eqref{ppg} satisfies the condition
\begin{eqnarray}\label{pp2}
P_{\xi,x^*}(x/\ox):=u_{x^*}(x)-\xi\|x-\ox\|\le P_{\xi,x^*}(\ox/\ox)=\la x^*,\ox\ra\;\mbox{ for all }\;x\in V(\ox)\cap\O
\end{eqnarray}
with $\xi:=\ve+\nu>\ve$. The latter tells us by definition \eqref{trap1} that $\ox$ is a local stationary trap relative to $\O$ in behavioral dynamics with the linear utility function $u_{x^*}(x)$ whenever $\xi>\ve$ for the weight factor $\xi$ in the proximal payoff \eqref{ppg}. $\h$\vspace*{0.05in}

The next theorem makes a connection between local stationary traps and the {\em extremal principle} for closed set systems that is yet another fundamental result of variational analysis with numerous applications; see, e.g.,  \cite{m06,m18} and the references therein. We primarily use here the (approximate) $\ve$-{\em extremal principle} from \cite[Definition~2.5(i)]{m06}, which holds in any Asplund space by \cite[Theorem~2.20]{m06}.

First we recall the notion of local extremal points of extremal systems of sets in normed spaces that is taken from \cite[Definition~2.1]{m06}.

\begin{Definition}{\bf (set extremality).}\label{ext} We say that $\ox$ is a {\sc locally extremal point} of the system of finitely many sets $\O_1,\ldots,\O_n$ in the normed space $X$ if:

{\bf (a)} $\ox$ is a common point of the sets $\O_i$ for all $i=1,\ldots,n$;

{\bf (b)} there exist a neighborhood $V$ of $\ox$ and sequences $\{a_{ik}\}\subset X$ for $i=1,\ldots,n$ such that $a_{ik}\to 0$ as $k\to\infty$ whenever $i=1,\ldots,n$ and that
\begin{eqnarray*}
\bigcap_{i-1}^n\Big(\O_i-a_{ik}\Big)\cap V=\emp\;\mbox{ for all large natural numbers }\;k.
\end{eqnarray*}
In this case we say that $\{\O_1,\ldots,\O_n,\ox\}$ is an {\sc extremal system} in $X$.
\end{Definition}

It has been well recognized in variational analysis and its applications that the concept of set extremality encompasses various notions of optimal solutions in problems of scalar, vector, and set-valued optimization, equilibria, game theory, systems control, models of welfare economics, etc. On the other hand, local extremal points naturally appear in developing generalized differential calculus, geometric aspects of functional analysis, and related disciplines. We refer the reader to both volumes of \cite{m06} and to the new book \cite{m18} for the extensive theoretical material on these issues and numerous applications.

Now we are ready to establish the aforementioned relationship between set extremality and stationary traps via the $\ve$-{\em extremal principle} in Asplund spaces.

\begin{Theorem}{\bf (local extremal points and stationary traps via the $\ve$-extremal principle).}\label{ext-trap} Let $\{\O_1,\ldots,\O_n,\ox\}$ be an extremal system in an Asplund space $X$, where the sets $\O_i$ are locally closed around $\ox$. Then for any $\ve>0$ there exist points
\begin{eqnarray}\label{x-cl}
x_i\in\O_i\;\mbox{ with }\;\|x_i-\ox\|\le\ve\;\mbox{ as }\;i=1,\ldots,n
\end{eqnarray}
and dual elements $x^*_i\in X^*$, which are rates of change for the linear utility functions
\begin{eqnarray}\label{ut}
u_{x^*_i}(x):=\la x^*_i,x\ra,\quad i=1,\ldots,n,
\end{eqnarray}
while satisfying the conditions
\begin{eqnarray}\label{ext1}
x^*_1+\ldots+x^*_n=0\;\mbox{ and }\;\|x^*_1\|+\ldots+\|x^*_n\|=1,
\end{eqnarray}
such that each $x_i$ is a local stationary trap relative to $\O_i$ with respect to the linear utility function $u_{x^*_i}(x)$ in \eqref{ut} for any weight factor $\xi>\ve$.
\end{Theorem}
{\bf Proof.} Since $\ox$ is a locally extremal point of the system $\{\O_1,\ldots,\O_n\}$ of locally closed sets $\O_i$ in the Asplund space $X$, we can employ the $\ve$-extremal principle from \cite[Theorem~2.20]{m06} and find $x_i$ satisfying \eqref{x-cl} and $x^*_i$ satisfying \eqref{ext1} such that
\begin{eqnarray}\label{e-nor}
x^*_i\in\Hat N_\ve(x_i;\O_i)\;\mbox{ for all }\;i=1,\ldots,n.
\end{eqnarray}
Then we use the stationary trap description from Proposition~\ref{nor-trap} of $\ve$-normals $x^*_i$ to $\O_i$ at points $x_i$ that are sufficiently close to the given local extremal point $\ox$. $\h$

\begin{Remark}{\bf (further discussions on geometric and analytic evaluations of stationary traps).}\label{disc} {\rm It is worth mentioning the following interpretations and extensions of the above evaluations of local stationary traps:

{\bf (i)} Taking into account the definitions of the advantage of change $A(\ox/x)$ from $\ox$ to $x$ in \eqref{A} and the corresponding cost of changing $C(\ox,x)$ in \eqref{C}, we can interpret the local stationary trap description of $\ve$-normals $x^*\in\Hat N_\ve(\ox;\O)$ from Proposition~\ref{nor-trap} as
\begin{eqnarray}\label{pp3}
A(\ox/x)\le\xi C(\ox,x)\;\mbox{ for all }\;v\in V(\ox)\cap\O
\end{eqnarray}
with the corresponding weight factor $\xi$ exceeding $\ve$, the linear utility/evaluation function $u_{x^*}(x)$ from \eqref{ut} with the rate of change $x^*$, and some neighborhood $V(\ox)$ of $\ox$ relative to $\O$ that depends on the factors above. Then \eqref{pp3} says that it is {\em not worthwhile to move} from the local stationary trap position $\ox$ to a point $x\in\O$ nearby.

{\bf (ii)} Given a payoff function $\ph\colon X\to\oR$, consider its {\em tilt perturbation}
\begin{eqnarray}\label{tilt}
\ph_v(x):=\ph(x)-\la v,x\ra\;\mbox{ with some }\;v\in X^*.
\end{eqnarray}
The importance of tilt perturbations and the related notion of {\em tilt stability}, introduced by Poliquin and Rockafellar \cite{pol-roc98}, have been well recognized in variational analysis and optimization; see, e.g., \cite{dru-lew13,mor-nghia15,mr12} with the references therein for more recent publications. In behavioral economics, tilt perturbations were proposed and developed by Thaler \cite{th08} in the framework of {\em acquisition utility}. Then the evaluation results obtained above in terms of $\ve$-subgradients \eqref{e-sub} can be interpreted in the way that, instead of minimizing the original payoff $\ph$, we actually minimize its {\em tilt-perturbed} counterpart \eqref{tilt} shifted by the {\em cost of changing/resistance} term $\xi\|x-\ox\|$. This reflects behavioral aspects well understood in psychology, where agents balance between {\em desirability} issues (minimizing $\ph_v(\cdot)$ in our case) and {\em feasibility} ones, i.e., minimizing their costs of changing $C(\ox,x)$.

{\bf (iii)} In various behavioral situations described by models in finite-dimensional or Hilbert spaces $(X^*=X)$ it is reasonable to replace the term $\|x-\ox\|$ in the costs of changing by its {\em quadratic} modification of the resistance to change $R=D[I]=C(\ox,x)^2$. Such settings can be investigated similarly to the above developments with replacing the collections of $\ve$-subgradients \eqref{e-sub} by the {\em proximal subdifferential}
\begin{eqnarray*}
\partial_P\vt(\ox):=\big\{v\in X\big|\;\vt(x)\ge\vt(\ox)+\la v,x-\ox\ra-(\rho/2)\|x-\ox\|^2\;\mbox{ for all }\;x\;\mbox{ near }\;\ox\big\}.
\end{eqnarray*}
This construction is useful for the design and investigation of {\em proximal-type algorithms} in optimization and related areas with applications to behavioral sciences; see, e.g., \cite{as11}.}
\end{Remark}

\section{Summary of Major Finding and Future Research}\sce

This paper reveals two-sided relationships between some basic notions and results of variational analysis with variational rationality in behavioral sciences. From one hand, we apply well-recognized constructions and principles of variational analysis and generalized differentiation to the study of stationary traps and related aspects of human dynamics. On the other hand, our new results provide valuable behavioral interpretations of general notions of variational analysis that is done for the first time in the literature.\vspace*{0.05in}

Among the {\em major finding} in this paper we underline the following:

$\bullet$ Introducing the notions of optimistic and pessimistic evaluations of payoff functions and derive efficient linear optimistic evaluation of the original payoffs and proximal payoffs via subgradients of convex analysis for original payoffs and $\ve$-subgradients of variational analysis for reference-dependent proximal payoff functions.

$\bullet$ Establishing a relationship between subgradient and weight factors that ensures simultaneously a linear evaluation of rates of change for proximal payoffs and a behavioral interpretation of $\ve$-subgradients for general extended-real-valued functions.

$\bullet$ Deriving certificates (sufficient conditions) for exact stationary traps in behavioral dynamics expressed in terms of $\ve$-subgradients of payoff functions.

$\bullet$ Similar subgradient certificates are obtained for approximate stationary traps.

$\bullet$ The obtained results allow us to determine and classify various special types of exact and approximate stationary traps in behavioral dynamics.

$\bullet$ By using powerful variational principles of Ekeland and lower subdifferential types, we establish the existence of approximate traps that satisfy certain optimality and subdifferential properties. This sheds light on the very nature of such traps under perturbations and on rates of change in behavioral dynamics at and around such positions.

$\bullet$ Besides the aforementioned analytic evaluations of exact and approximate stationary traps by using subgradients, we develop their geometric evaluations in terms of generalized normals to closed sets. These developments allow us also to describe generalized normals of variational analysis via stationary traps in behavioral dynamics with linear utility functions and high enough weight factors in costs of changing.

$\bullet$ Finally, we relate local stationary traps in behavioral processes described via linear utility functions with local extremal points of set systems and the fundamental extremal principle of variational analysis that plays a crucial role in both theory and applications.\vspace*{0.05in}

The established two-sided relationships between variational analysis and variational rationality in behavioral sciences open the gate for further developments in this direction, which we plan to pursue in our {\em future research}. Besides the questions mentioned in the remarks above, they include while are not limited to:

$\bullet$  Considering {\em variational traps} in behavioral dynamics, which indicate positions that are worthwhile to approach and reach by a successions of worthwhile moves, but not worthwhile to quit. Helpful insights into the study of variational traps are given by the constructive dynamic {\em proofs} of the extremal principle in both finite and infinite dimensions; see \cite[Theorems~2.8 and 2.10]{m06}.

$\bullet$ Implementing the established relationships and proposed variational ideas into developing {\em numerical algorithms} of proximal, subgradient, and majorization-minimization type with applications to behavioral science modeling; compare, e.g., \cite{as11,bp16}.

$\bullet$ Developments and applications of the extremal principle and related tools of variational analysis to problems with {\em variable cone preferences}, which naturally arise in behavioral sciences via the variational rationality approach; see \cite{bms15a,bms15b} and also compare it with the books \cite{aky18,eic14,ktz15} treated general multiobjective optimization problems of this type.


\begin{thebibliography}{99}

\bibitem{aky18} Q. H. Ansari, E. K\"obis and J.-C. Yao, {\em Vector Variational Inequalities and Vector Optimization. Theory and Applications}, Springer,
Berlin, 2018.

\bibitem{as11} H. Attouch and A. Soubeyran, Local search proximal algorithms as decision dynamics with costs to move, {\em Set-Valued Var. Anal.} {\bf 19}
(2011), 157--177.

\bibitem{bao-mor10} T. Q. Bao and B. S. Mordukhovich, Set-valued optimization in welfare economics, {\em Adv. Math. Econ.} {\bf 13} (2010), 114--153.

\bibitem{bms15a} T. Q. Bao, B. S. Mordukhovich and A. Soubeyran, Variational analysis in psychological modeling, {\em J. Optim. Theory Appl.}
{\bf 164} (2015), 290--315.

\bibitem{bms15b} T. Q. Bao, B. S. Mordukhovich and A. Soubeyran, Fixed points and variational principles with applications to capability theory of wellbeing
via variational rationality, {\em Set-Valued Var. Anal.} {\bf 23} (2015), 375--398.

\bibitem{bp16} J. Bolte and E. Pauwels, Majoration-minimization procedures and convergence of SQP methods for semi-algebraic and tame programs,
{\em Math. Oper. Res.} {\bf 41} (2016), 442--465.

\bibitem{dru-lew13} D. Drusvyatskiy and A. S. Lewis, Tilt stability, uniform quadratic growth, and strong metric regularity of
the subdifferential, {\em SIAM J. Optim.} {\bf 23} (2013), 256--267.

\bibitem{eic14} G. Eichfelder, {\em Variable Ordering Structures in Vector Optimization}, Springer, Berlin, 2014.

\bibitem{kt79} D. Kahneman and A. Tversky, Prospect theory: an analysis of decision under risk, {\em Economet.} {\bf 47} (1979), 263--291.

\bibitem{ktz15} A. A. Khan, C. Tammer and C. Z\u{a}linescu, {\em Set-Valued Optimization. An Introduction with Applications}, Springer, Berlin. 2015.

\bibitem{m06} B. S. Mordukhovich, {\em Variational Analysis and Generalized Differentiation, I: Basic Theory; II: Applications}, Springer, Berlin, 2006.

\bibitem{m18}  B. S. Mordukhovich, {\em Basics of Variational Analysis and Applications}, Springer, New York, 2018.

\bibitem{mor-nghia15} B. S. Mordukhovich and T. T. A. Nghia, Second-order characterizations of tilt stability with applications to nonlinear programming,
{\em Math. Program.} {\bf 149} (2015), 83--104.

\bibitem{mr12} B. S. Mordukhovich and R. T. Rockafellar, Second-order subdifferential calculus with applications to tilt stability in optimization, {\em
SIAM J. Optim.} {\bf 22} (2012), 953--986.

\bibitem{pol-roc98} R. A. Poliquin and R. T. Rockafellar, Tilt stability of a local minimum, {\em SIAM J. Optim.} {\bf 8} (1998), 287--299.

\bibitem{rw} R. T. Rockafellar and R. J-B. Wets, {\em Variational Analysis}, Springer, Berlin, 1998.

\bibitem{s09} A. Soubeyran, {\em Variational Rationality: A Theory of Individual Stability and Change, Worthwhile and Ambidextry Behaviors},
GREQAM-AMSE, Aix-Marseille University, 2009.

\bibitem{s10} A. Soubeyran, {\em Variational Rationality and the Unsatisfied Man: Routines and the Course Pursuit between Aspirations, Capabilities and
Beliefs}, GREQAM-AMSE, Aix-Marseille University, 2010.

\bibitem{s18a} A. Soubeyran, {\em Variational Rationality: Proximal Dynamics and Stationary Traps}, GREQAM-AMSE, Aix-Marseille  University, 2018.

\bibitem{s18b} A. Soubeyran, {\em Variational Rationality: The Formation of Preferences and Intentions}, GREQAM-AMSE, Aix-Marseille University, 2018.

\bibitem{th08} R. H. Thaler, Mental accounting and consumer choice, {\em Marketing Sci.} {\bf 27} (2008), 15--25.

\bibitem{tk91} A. Tversky and D. Kahneman, Loss aversion in riskless choice: a reference-dependent model, {\em Quarterly J. Econ.} {\bf 106} (1991),
1039--1061.




\end{thebibliography}
\end{document}